\documentclass{article}[12pt]

\usepackage{amsmath}
\usepackage{amssymb}
\usepackage{amsthm}
\usepackage{amsfonts}
\usepackage{graphicx}
\usepackage{mathrsfs}
\usepackage{enumerate}

\newtheorem{thm}{Theorem}
\newtheorem{lem}{Lemma}
\newtheorem{fac}{Fact}
\newtheorem{cor}{Corollary}

\newtheorem{defn}{Definition}
\newtheorem{rem}{Remark}

\newenvironment{que}{\vskip2mm \noindent{\bf Question}\ }{\vskip2mm}

\newcommand{\im}{invariant measure}
\newcommand{\sq}{sequence} 

\newcommand{\z}{\mathbb Z}
\newcommand{\na}{\mathbb N}

\newcommand{\p}{\mathscr P}              
\newcommand{\Q}{\mathscr Q}              
\newcommand{\R}{\mathscr R}              
\newcommand{\xsmt}{$(X,\mathfrak B,\mu,T)$}
\newcommand{\xt}{$(X,T)$}
\newcommand{\ys}{$(Y,S)$}
\newcommand{\xsm}{$(X,\mathfrak B,\mu)$}
\newcommand{\ycns}{$(Y,\mathfrak C,\nu,S)$}

\newcommand{\tl}{topological}
\newcommand{\ds}{dynamical system}

\title{Measure-theoretic chaos\footnotetext{Research of the first author supported from resources 
for science in years 2009-2012 as research project (grant MENII N N201 394537, Poland) \newline
\it Mathematics Subject Classification (2010): \rm 37A35, 37B40 \newline
\it Keywords and phrases: \rm measure-theoretic chaos, distributional chaos, ergodic system, scrambled set, positive entropy.}
}

\date{}

\author{Tomasz Downarowicz and Yves Lacroix}
\begin{document}
\maketitle

\begin{abstract}
We define new isomorphism-invariants for ergodic measure-preserving systems on standard probability spaces, called measure-theoretic
chaos and measure-theoretic$^+$ chaos. These notions are analogs of the topological chaoses {\rm DC2} and 
its slightly stronger version (which we denote by {\rm DC}{\small$1\!\tfrac12$}). We prove that: 1. If a \tl\ system is measure-theoretically (measure-theoretically$^+$) chaotic with respect to at least one of its ergodic 
measures then it is \tl ly {\rm DC2} ({\rm DC}{\small$1\!\tfrac12$}) chaotic. 2.~Every ergodic system with positive Kolmogorov--Sinai entropy is measure-theoretically$^+$ chaotic (even in a bit stronger uniform sense). We provide an example showing that the latter statement cannot be reversed, i.e., of 
a system of entropy zero with uniform measure-theoretic$^+$ chaos.

\end{abstract}

\bigskip
\section{Introduction} 
The notion of chaos was invented by Li and Yorke in their seminal paper \cite{LY} in the context of continuous transformations of the interval. Since then several refinements of chaos have been introduced and extensively studied, for instance three versions of so-called \emph{distributional chaos} ({\rm DC1}, {\rm DC2} and {\rm DC3}) invented by Sm\'ital et al (\cite{SS}, \cite{SSt}, \cite{BSS}). All these notions refer to topological dynamical systems (actions of the iterates of a single continuous transformation $T$ on a compact metric space $X$) and strongly rely on the observation of distances between orbits, and the existence of so-called \emph{scrambled pairs} (or \emph{scrambled sets} -- usually uncountable). There are other notions of chaos, such as Devaney chaos or omega chaos, defined without the notion of scrambling -- these are not addressed in our paper.

Unlike in the case of most other notions in dynamics, there have been, to our knowledge, no 
successful attempts to create a measure-theoretic analog of chaos -- a notion applicable to measure-preserving transformations of a standard probability space without any specified topology. 
Although ``measure-theoretic chaos'' appears in the titles of some papers (e.g.\,\cite{WuWang}),
it still applies to topological systems. There are two major reasons why, at a first glance, 
it seems difficult to create such an analog: 
\begin{itemize}
	\item A standard probability space can be modeled as a compact metric space in many different ways. A pair (or set) scrambled in one metric need not be scrambled in another. 
\end{itemize}
\begin{itemize}	
	\item A scrambled set in a topological dynamical system very often has measure zero for every \im. It is always so for example in case of distributional chaos (in any version) -- we will explain this later. Li--Yorke-scrambled sets can have positive measure or even be equal to the entire space, but systems with such large scrambled sets are rather exceptional, and any notion of measure-theoretic chaos based on the analogy to these systems would be very restrictive (comp. \cite{WuWang}).	In all other cases, a scrambled set can be easily added to the space (or discarded from it) in a way negligible from the point of view of measure. In other words, chaos based on the existence of a scrambled set is not stable under measure-theoretic isomorphisms. 
\end{itemize}

Inspired by the methods developed in \cite{D}, in this note we propose a way to overcome these difficulties. We define chaos in measure-theoretic systems using exclusively the measurable structure of the space, and so that it becomes an invariant of measure-theoretic isomorphism. Our new notions maintain their original character -- they are defined in terms of uncountable scrambled sets. Moreover, they are related to their topological prototypes and also to positive entropy exactly as one would expect (we will give more details in a moment). 

Among the topological notions of chaos we have chosen one -- the distributional chaos {\rm DC2} (with variants) -- as the starting point to define its measure-theoretic analog. This new notion, which we call simply \emph{the measure-theoretic chaos}\footnote{
We have decided to suppress the adjective ``distributional'' because ``distribution'' is a synonym of ``measure''. In case of \tl\ chaos, this adjective indicates ``some reference to measures'' (maintaining reference to the metric), while here we have reference to a measure and nothing else in fact, so the adjective ``measure-theoretic'' should suffice.
}, 
meets all our expectations regarding its relations with the topological prototype, and it inherits the most important implications 
between chaos and entropy.

Other notions of chaos are not so well adaptable to the measure-theoretic context; the attempted analogs of Li--Yorke and DC1 chaoses
fail a key property allowing to prove that they imply their \tl\ prototypes (see Remark \ref{fail} for {\rm DC1}). It is possible to copy our scheme for {\rm DC3} (see Remark~\ref{dc3}), but because generally this notion is very weak (it can occur even in distal systems), we have decided to skip it. Nonetheless, for completeness of the survey in the next section, we include the definitions of Li--Yorke, {\rm DC1} and {\rm DC3} chaoses in \tl\ systems.

Let us recall that Blanchard, Glasner, Kolyada and Maass have proved that positive topological entropy implies Li--Yorke chaos (see \cite{BGKM}). This result has been recently strengthened by the first author of this note: positive topological entropy implies distributional chaos {\rm DC2} (see \cite{D}). Let us also recall that for interval maps all three versions of distributional chaos ({\rm DC1}, {\rm DC2} and {\rm DC3}) are equivalent to positive topological entropy (see \cite{SS}). We can now be more specific about maintaining these implications by our new notion. We will prove that:

\begin{itemize}
	\item A \tl\ \ds\ which is measure-theoretically chaotic with respect to at least one of its \im s is {\rm DC2} chaotic.
	\item A measure-theoretic system with positive Kolmogorov--Sinai entropy is measure-theoretically chaotic
	(in particular, a topological system with positive topological entropy is measure-theoretically chaotic for at least one of its \im s --  and thus {\rm DC2} chaotic). 
	\item For a continuous transformation of the interval, positive topological entropy is equivalent to measure-theoretic chaos for some of its \im s.
\end{itemize}
The last statement is a direct consequence of the preceding two, the variational principle and the equivalence between {\rm DC2} and positive topological entropy for interval maps, thus we do not need prove it separately. We believe that the above assembly of relations (plus the fact that our notion is an isomorphism invariant) is a good enough reason to consider our notion a successful analog of distributional chaos in measure-theoretic dynamics.

\section{Review of topological chaos}

Let us begin with a review of topological notions of chaos: Li--Yorke, {\rm DC1}, {\rm DC2} and {\rm DC3}. Later we will also introduce a notion intermediate between {\rm DC1} and {\rm DC2}, which we denote {\rm DC}{\small$1\!\tfrac12$}. All these notions are defined in the same manner: \emph{there exists an uncountable scrambled set}, where a scrambled set is one whose every pair of distinct elements is scrambled. The only remaining detail is the meaning of a ``scrambled pair'' for the above types of chaos. The definitions given below are equivalent to the most commonly appearing in the literature but expressed using a slightly different language (this change is meant for an easy adaptation to the measure-theoretic situation). We will also define uniform versions of {\rm DC1} and {\rm DC2} (and later -- of {\rm DC}{\small$1\!\tfrac12$}).

\medskip
Traditionally, a pair $(x,y)$ is Li--Yorke-scrambled if 
$$
\liminf_{n\to\infty} d(T^nx,T^ny)=0\text{ \ \ and \ \ }\limsup_{n\to\infty} d(T^nx,T^ny)>0.
$$ 
This can be rephrased as follows:
\begin{itemize}
	\item A pair $(x,y)$ is \emph{Li--Yorke-scrambled} if there exist: an increasing \sq\ $n_i$ such that $d(T^{n_i}x,T^{n_i}y)\overset{i}\longrightarrow0$, 
another increasing \sq\ $m_i$, and a positive number $s>0$, such that $d(T^{m_i}x,T^{m_i}y)\ge s$ for all $i$.
\end{itemize}
Distributional scrambling {\rm DC1} and {\rm DC2} are similar, except we put density constraints 
on the sequences $n_i$ and $m_i$:
\begin{itemize}
	\item A pair $(x,y)$ is \emph{{\rm DC1}-scrambled} if there exist: an increasing \sq\ $n_i$ of upper density 1, such that $d(T^{n_i}x,T^{n_i}y)\overset{i}\longrightarrow0$, 
	another increasing \sq\ $m_i$ of upper density 1, and a number $s>0$, such that $d(T^{m_i}x,T^{m_i}y)\ge s$ for all $i$.
  \item A pair $(x,y)$ is \emph{{\rm DC2}-scrambled} if there exist: an increasing \sq\ $n_i$ of upper density 1, such that $d(T^{n_i}x,T^{n_i}y)\overset{i}\longrightarrow0$, another increasing \sq\ $m_i$ of \emph{positive} upper density, and a number $s>0$, such that $d(T^{m_i}x,T^{m_i}y)\ge s$ for all $i$.
\end{itemize}
The resulting chaos {\rm DC1} is called \emph{uniform} if the constant $s$ can be chosen common for all pairs in the scrambled set. For uniformity of {\rm DC2} we will require that all pairs in the scrambled set have common both the parameter $s$ and a positive lower bound $\eta$ for the upper density of the \sq s $m_i$. 

Scrambling for {\rm DC3} has a slightly different structure:
\begin{itemize}  
  \item A pair $(x,y)$ is \emph{{\rm DC3}-scrambled} if there exists an interval $(a,b)$ such that for every $s\in (a,b)$ 
  the \sq\ of the times $n$ when $d(T^{n}x,T^{n}y)\ge s$ does not have a density (upper and lower densities differ).\footnote{
Traditionally, distributional scrambling is defined using the functions 
$\Phi^*_{x,y}(t)$ and $\Phi_{x,y}(t)$ defined for $t>0$ as, respectively, the upper and lower densities of the set of times $n$ when $d(T^nx,T^ny)< t$. Clearly, $\Phi^*_{x,y}\ge \Phi_{x,y}$, both functions increase with $t$, reaching the value $1$ for the diameter of $X$. One can define both functions at zero as the limit values as $t\to 0^+$. A pair $(x,y)$ is  
\begin{itemize}
\item {\rm DC1}-scrambled if $\Phi_{x,y}^*(0)=1$ and $\Phi_{x,y}(s)=0$ for some $s$; 
\item {\rm DC2}-scrambled if $\Phi_{x,y}^*(0)=1, \Phi_{x,y}(0)<1$;
\item {\rm DC3}-scrambled if $\Phi_{x,y}^*(s)>\Phi_{x,y}(s)$ on an open interval. 
\end{itemize}
}
\end{itemize}

Requesting the upper density of $m_i$ to be arbitrarily close to 1 we produce a notion intermediate between {\rm DC2} and {\rm DC1}:
\begin{itemize} 
\item A pair $(x,y)$ is \emph{{\rm DC}{\small$1\!\tfrac12$}-scrambled} if there exist: an increasing \sq\ $n_i$ of upper density 1, such that $d(T^{n_i}x,T^{n_i}y)\overset{i}\longrightarrow0$, 
and, for every $\eta<1$, an increasing \sq\ $m_{\eta,i}$ of upper density at least $\eta$, and a number $s_\eta>0$, such that $d(T^{m_{\eta,i}}x,T^{m_{\eta,i}}y)\ge s_\eta$ for all $i$.\footnote{
In other words, a pair $(x,y)$ is {\rm DC}{\small$1\!\tfrac12$}-scrambled if $\Phi_{x,y}^*(0)=1, \Phi_{x,y}(0)=0$. 
}
\end{itemize}
The meaning of \emph{chaos {\rm DC}{\small$1\!\tfrac12$}} is clear: there exists an uncountable {\rm DC}{\small$1\!\tfrac12$}-scram\-bled set. \emph{Uniform chaos {\rm DC}{\small$1\!\tfrac12$}} occurs when the function $\eta\mapsto s_\eta$ is common for all pairs in the scrambled set\footnote{
In \cite{D} it is proved that positive topological entropy implies uniform chaos {\rm DC2}. In this paper we will strengthen that result: positive topological entropy implies uniform chaos {\rm DC}{\small$1\!\tfrac12$}. This is why we think uniform {\rm DC}{\small$1\!\tfrac12$} is worth a separate formulation. 
Similar result is obtained in this paper for the measure-theoretic analog.
}.
It is easy to see that 
\begin{gather*}
\text{{\rm DC1}} \implies \text{{\rm DC}{\small$1\!\tfrac12$}} \implies\text{{\rm DC2}} \implies\text{{\rm DC3}},\\ 
\text{{\rm DC2}}\implies\text{Li--Yorke chaos,}
\end{gather*} 
and that {\rm DC1} through {\rm DC2} (including the uniform versions) are invariants of \tl\ conjugacy (see \cite{SSt}).
{\rm DC3} does not imply Li--Yorke and, as we mentioned earlier, is not a \tl\ invariant. 
There are easy examples showing that {\rm DC1} is essentially stronger than {\rm DC}{\small$1\!\tfrac12$}, in fact
(as we will show later) every system with positive topological entropy is uniformly {\rm DC}{\small$1\!\tfrac12$}, while Piku\l a
provided an example of a system with positive topological entropy which is not {\rm DC1} (\cite{P}).
It is not very hard to construct an example showing that {\rm DC2} (even uniform) is essentially 
weaker than {\rm DC}{\small$1\!\tfrac12$}. We refrain from providing such an example in this note. A more 
interesting question is whether {\rm DC2} (or uniform {\rm DC2}) 
\emph{persistent under removing null sets} (see the formulation of Theorem~\ref{persistent}) 
implies {\rm DC}{\small$1\!\tfrac12$}. At the moment we leave this problem open, with a conjecture that the answer is negative.


\medskip
Notice that the condition for a {\rm DC2}-scrambled pair has a beautiful translation to the language 
of ergodic averages. A pair $(x,y)$ is {\rm DC2}-scrambled if and only if 
$$
\liminf_{n\to\infty} \frac1n \sum_{i=1}^n d(T^ix,T^iy)=0 \ \ \text{\ \ and \ } \ \ \limsup_{n\to\infty} \frac1n \sum_{i=1}^n d(T^ix,T^iy)>0.
$$
Note how this formulation is analogous to the original condition for Li--Yorke scrambling. It is not hard to see that uniformity of {\rm DC2} is equivalent to the upper limit seen above on the right having a common positive lower bound for all pairs in the scrambled set. The fact that {\rm DC2} can be phrased in terms of ergodic averages makes it the best and most natural candidate to become a base for creating a measure-theoretic analog. 

\section{Definitions of measure-theoretic chaos}
At this point we leave the topological setup of a compact metric space and we move into the context of a standard probability space \xsm, where $\mathfrak B$ is a complete sigma-algebra and $\mu$ is a probability measure on $\mathfrak B$, on which we consider the action of a measure-preserving transformation $T$. As it was said before, in order to define measure-theoretic chaos we must overcome two difficulties, the first of which is that the definition of scrambling must not refer to any metric. This is done using refining sequences of finite partitions.

\begin{defn}
A \sq\ of finite measurable partitions $(\mathscr P_k)_{k\ge 1}$ is called \emph{refining} if 
$\mathscr P_{k+1}\succcurlyeq \mathscr P_k$ for every $k$ and jointly they generate $\mathfrak B$
(i.e., $\mathfrak B$ is the smallest complete sigma-algebra containing all the partitions $\mathscr P_k$).
\end{defn}

\begin{defn}
Fix a refining \sq\ of finite measurable partitions $(\mathscr P_k)$. A pair of points $(x,y)$ 
is \emph{$(\mathscr P_k)$-scrambled} if
\begin{itemize}
	\item There exists a sequence $n_i$ of upper density 1 such that for every $k$ and large enough $i$,
	$T^{n_i}x$ belongs to the same atom of $\mathscr P_k$ as $T^{n_i}y$.
	\item There exists a sequence $m_i$ of positive upper density, and $k_0$ such that, for every $i$,
	$T^{m_i}x$ and $T^{m_i}y$ belong to different atoms of $\mathscr P_{k_0}$.
\end{itemize}
\end{defn}

The second major difficulty is to assure that our chaos is an isomorphism invariant. This is achieved
by requiring the existence of a scrambled set for \emph{every} refining \sq\ of finite partitions, as it is done in the definitions given below. We will show in the next section that the notions of chaos so constructed are indeed isomorphism invariants.

\begin{defn} \label{ch}
A measure-preserving transformation $T$ of a standard probability space \xsm\ is \emph{measure-theoretically  
chaotic} if for every refining \sq\ of finite partitions $(\mathscr P_k)$ there exists an uncountable $(\mathscr P_k)$-scrambled 
set.
\end{defn}

\begin{defn} \label{uch}
The above defined chaos is \emph{uniform} if (for any refining \sq\ $(\p_k)$) all distinct pairs in the scrambled set are 
$(\mathscr P_k)$-scrambled with a common parameter $k_0$ and with a common positive lower bound $\eta$ on the upper density of 
the \sq s $m_i$.
\end{defn}

It is also easy to define a stronger version of measure-theoretic chaos, an analog of {\rm DC}{\small$1\!\tfrac12$}. Precisely this 
type of chaos is implied by positive entropy, hence we find it worth a presentation. It suffices to modify the definition of 
scrambled pairs:
\begin{defn} A pair $(x,y)$ is $(\mathscr P_k)^+$-scrambled if
\begin{itemize}
	\item There exists a sequence $n_i$ of upper density 1 such that for every $k$ and large enough $i$,
	$T^{n_i}x$ belongs to the same atom of $\mathscr P_k$ as $T^{n_i}y$.
	\item For every $\eta>0$ there exists a sequence $m_{\eta,i}$ of upper density at least
	$\eta$, and $k_\eta$ such that, for every $i$, $T^{m_{\eta,i}}x$ and $T^{m_{\eta,i}}y$ belong to different atoms of $\mathscr P_{k_\eta}$.
\end{itemize}
\end{defn}
Replacing $(\mathscr P_k)$-scrambling in the definition of the measure-theoretic chaos by 
$(\mathscr P_k)^+$-scrambling we obtain \emph{measure-theoretic$^+$ chaos}. For \emph{uniform measure-theoretic$^+$ chaos} we require that the function $\eta\mapsto k_\eta$ is common for all pairs in the scrambled set.
\bigskip

\section{Presentation of the measure-theoretic chaos}

In this section we formulate and prove our statements concerning the notion of measure-theoretic chaos, in particular
its isomorphism invariance and its relations with the topological counterpart. 

\begin{thm}\label{inva}
Suppose the systems \xsmt\ and \ycns\ are isomorphic. Then \xsmt\ is measure-theoretically  
(uniformly measure-theoretically,\break measure-theoretically$^+$, uniformly measure-theoretically$^+$) chaotic if and only 
if so is \ycns.
\end{thm}

\begin{proof} Let $\pi:X\to Y$ be the isomorphism. Recall that $\pi$ is a measurable bijection between full sets $X'\subset X$ 
and $Y'\subset Y$ (i.e., sets of full measure in the respective spaces), intertwining the actions of $T$ and $S$. By a standard argument we can arrange that $X'$ and $Y'$ are invariant, that is $T(X')\subset X'$ and $S(Y')\subset Y'$. Let $(\Q_k)$ be an arbitrarily chosen refining \sq\ of partitions of $Y$. Denote by $\Q'_k$ the restriction of $\Q_k$ to $Y'$ and let $\p'_k$ be the partition of $X'$ obtained as the preimage by $\pi$ of $\Q'_k$. Finally let $\p_k$ denote the partition of $X$ consisting of the elements of $\p'_k$ and the null set $C=X\setminus X'$. It is obvious that $(\p_k)$ is a refining \sq\ of partitions in $X$. If the system on $X$ is chaotic (in any of the four considered senses) then there exists an uncountable $(\p_k)$-scrambled set $E$ (for the corresponding meaning of scrambling). By $(\p_k)$-scrambling and invariance of $X'$, for every $x\in E$ there exists $n_x$ such 
that $T^{n_x}x\in X'$. Since $E$ is uncountable, it has an uncountable subset $E'$ with a common $n_x$. 
The set $E''=T^{n_x}(E')$ is uncountable, $(\p'_k)$-scrambled and contained in $X'$. Now the set $\pi(E'')$ is obviously $(\Q_k')$-scrambled in $Y'$ which immediately implies that it is $(\Q_k)$-scrambled in $Y$. This ends the proof.
\end{proof}

Although our notions of chaos formally apply to all measure-theoretic systems, we will focus on the most important, ergodic, case. Most of the theorems stated below require ergodicity anyway. Thus, throughout the remainder of the paper we will assume that $\mu$ is ergodic.

In this context we will provide conditions equivalent to measure-theoretic chaos (and its variants) referring to only 
one refining \sq\ of partitions and ``persistence under removing null sets'':

\begin{thm}\label{nullset} Let $(\p_k)$ be a fixed refining \sq\ of finite partitions of $X$. The ergodic system \xsmt\ is measure-theoretically (measure-theoretically$^+$, uniformly measure-theoretically, uniformly measure-theoretically$^+$) chaotic 
if and only if for any null set $A$ (i.e., of measure zero) there exists an uncountable $(\p_k)$-scrambled ($(\p_k)^+$-scrambled, uniformly $(\p_k)$-scrambled, uniformly $(\p_k)^+$-scrambled) set disjoint from $A$.
\end{thm}

\begin{proof}
One implication is trivial, since removing a null set is in fact an isomorphism. We will focus on the nontrivial implication.
Suppose that no matter what null set is removed from $X$, there remains a $(\mathscr P_k)$-scrambled set. Consider another refining 
\sq\ of finite partitions $(\mathscr P'_{k'})$. Let $A$ be the null set such that all remaining points satisfy the assertion of the ergodic theorem with regard to all (countably many) elements of the field (attention, not sigma-field, just field) $\mathscr F$ generated by the partitions $(\mathscr P_k)$ and $(\mathscr P'_{k'})$ (that means, the orbit of every remaining point visits every 
set $B\in\mathscr F$ along a set of times whose density equals $\mu(B)$). By assumption, there exists an uncountable $(\mathscr P_k)$-scrambled set disjoint from $A$. We will show that the same set is $(\mathscr P'_{k'})$-scrambled. Take a pair of distinct points $(x,y)$ from this set and fix some $k'$. For an arbitrarily small $\delta>0$ there exists $k$ and a set $B$ of measure at most $\delta$ such that relatively on $X\setminus B$, the partition $\mathscr P'_{k'}$ is refined by $\mathscr P_k$. Clearly, $B$ belongs to $\mathscr F$. We know that the sequence of times $n$ when $T^nx$ and $T^ny$ belong to the same element of $\mathscr P_k$ has upper density 1. Removing the sequence of times when at least one of the above two points falls into $B$ we obtain a \sq\ of times of upper density at least $1-2\delta$, when the two points fall in the same element of $\mathscr P_k$ within $X\setminus B$ (and hence they fall into the same element of $\mathscr P'_{k'}$). Because $\delta$ was arbitrarily small, we get that $T^nx$ and $T^ny$ fall into the same element of $\mathscr P'_{k'}$ for times $n$ with upper density 1. Further, we know that there exists an index $k_0$ such that the sequence of times $n$ for which $T^nx$ and $T^ny$ belong to different elements of $\mathscr P_{k_0}$ has positive upper density, say $\eta$. Fix some $\delta>0$ much smaller than $\eta$ and find $k'_0$ and a set $C\in\mathscr F$ of measure at most $\delta$ such that $\mathscr P'_{k'_0}$ refines $\mathscr P_{k_0}$ relatively on $X\setminus C$. Removing the sequence of times when at least one of the above points falls into $C$, we obtain a \sq\ of upper density at least $\eta-2\delta$ (which is positive) when the two points fall into different elements of $\mathscr P_{k_0}$ within $X\setminus C$ (hence they are in different atoms of $\mathscr P'_{k'_0}$). 
We have proved that the pair $(x,y)$ is $(\mathscr P'_{k'})$-scrambled. 

For measure-theoretic$^+$ chaos it suffices to note that if $\eta$ is close to 1, so is $\eta-2\delta$. 

The same proof applies also to uniform measure-theoretic  chaos: if $k_0$ and $\eta$ are common to all pairs in the $(\mathscr P_k)$-scrambled set, the proof produces common parameters $k'_0$ and $\eta-2\delta$ for all pairs in the same set regarded as $(\mathscr P'_{k'})$-scrambled.

The argument for uniform measure-theoretic$^+$ chaos is the same as for the uniform 
measure-theoretic chaos, applied separately for every $\eta$ (with $k_\eta$ and $k'_\eta$ in place of $k_0$ and $k'_0$, respectively).
\end{proof}

\begin{rem}\label{fail}
Even if $\eta=1$, the proof produces $\eta-2\delta<1$. This is the reason why we gave up defining an analog of {\rm DC1}; 
it could not be tested using one \sq\ of partitions. Similarly, it would probably not imply {\rm DC1} in \tl\ systems (the proof of Theorem \ref{topol} as it is would not pass).
\end{rem}

The next theorem replaces the ``persistence under removing null sets'' by a much stronger property: ``ubiquitous presence of chaos'': scrambled sets exist inside any set of positive measure:

\begin{thm}\label{posset}
Let \xsmt\ be an ergodic measure-theoretically (uniformly measure-theoretically, measure-theoretically$^+$, uniformly measure-theoretically$^+$) chaotic system. Let $(\p_k)$ be a refining \sq\ of finite partitions and let $B\in\mathfrak B$ be a set 
of positive measure. Then there exists an uncountable $(\p_k)$-scrambled ($(\p_k)^+$-scrambled, uniformly $(\p_k)$-scrambled, 
uniformly $(\p_k)^+$-scrambled) set contained in $B$.
\end{thm}

\begin{proof}
Throughout the proof ``scrambled set'' stands for either $(\p_k)$-scrambled set, $(\p_k)^+$-scrambled set, uniformly $(\p_k)$-scrambled set, or uniformly $(\p_k)^+$-scrambled set, depending on the considered version of chaos.
Let $A\subset X$ be a measurable set containing no uncountable scrambled sets. We need to show that $\mu(A)=0$. Since the image of a scrambled set is obviously scrambled and has the same cardinality, $T^{-n}(A)$ does not contain uncountable scrambled sets either. This implies that $A' = \bigcup_{n=0}^\infty T^{-n}(A)$ does not contain uncountable scrambled set (otherwise an uncountable subset of the scrambled set would have to fall in one item of the union). But $A'$ is subinvariant (contains its preimage), hence, by ergodicity, its measure is either 1 or 0. The first possibility is excluded by Theorem \ref{nullset}. It follows that $\mu(A')=0$, in particular $\mu(A)=0$. 
\end{proof}

Let us devote one page to better understanding the phenomenon of ``ubiquitous presence of chaos''. This phenomenon is the major difference between how topological and measure-theoretic chaoses are constructed (making the latter much stronger). First of all, let us realize that such presence cannot be achieved by the existence of a scrambled set of full (or even positive) measure. Fact is, every $(\mathscr P_k)$-scrambled set must be a null set. The same applies to distributionally scrambled sets in topological \ds s: 

\begin{fac} Let $(\mathscr P_k)$ be a refining \sq\ of finite measurable partitions. Assume that $\mu$ is nonatomic. Then any $(\mathscr P_k)$-scrambled set has measure zero. Similarly, any {\rm DC3}-scrambled (and thus also {\rm DC2}-scrambled or {\rm DC1}-scrambled) set in a \tl\ \ds\ is a null set for all nonatomic \im s.
\end{fac}
\begin{proof}
Let $(x,y)$ be a $(\mathscr P_k)$-scrambled pair and let $k_0$ be the index in the definition of $(\mathscr P_k)$-scrambling. Consider the two-element partition of $X\times X$ into two sets: $\bigcup_{A\in\mathscr P_{k_0}}A\times A$ and its complement, $\bigcup_{A,B\in\mathscr P_{k_0}, A\neq B}A\times B$. The orbit of the pair visits the first 
set with upper density 1 and the other with positive upper density, so the visits in these sets do not have densities. Such pairs are exceptional (belong to a null set depending on the index $k_0$) for every ergodic measure on $X\times X$. Since there are countably many choices of $k_0$, the collection of all $(\mathscr P_k)$-scrambled pairs is a null set for any such measure, and hence also for any $T\times T$-\im, in particular, for $\mu\times\mu$. So, if $E$ is a $(\mathscr P_k)$-scrambled set, we have 
$(\mu\times\mu)(E\times E\setminus\Delta)=0$. Since $\mu$ is nonatomic, also $(\mu\times\mu)(\Delta)=0$, which implies
$(\mu\times\mu)(E\times E)=0$ and hence $\mu(E)=0$. 

The proof for {\rm DC3}-scrambled sets in \tl\ systems is identical, except that the two-set partition consists of the $s$-neighborhood of the diagonal and its complement.
\end{proof}

It is clear that the ``ubiquitous presence of chaos'' requires the union of all scrambled sets to be a set of full measure. Moreover, by a simple transfinite argument, there must exist a disjoint family of scrambled sets whose union is a full measure set. 
But even this last condition seems to be insufficient. Although we do not have an example of a dynamical system, it is easy to imagine an abstract family of disjoint uncountable null sets whose union has full measure, yet, this measure is supported by a set selecting only countably many points (or just one point) from each member of the family. Then by removing the rest (which is a null set) we destroy all the uncountable sets. So, the ``ubiquitous presence of chaos'' requires, most likely, an even more sophisticated configuration of the scrambled sets (than just the existence of a disjoint collection forming a full set). We give up further attempts to find an equivalent condition. What we have just learned for sure is that it is related to \emph{abundance} of scrambled sets rather than their individual largeness. 

\medskip
Next, we take care of the relations between the notions of measure-theoretic chaos and their topological prototypes. 

\begin{thm}\label{topol} Let $(X,T)$ be a \tl\ \ds\ and let $\mu$ be an ergodic $T$-\im. If the measure-theoretic system \xsmt\ 
(where $\mathfrak B$ denotes the Borel sigma-algebra completed with respect to $\mu$) is measure-theoretically (uniformly measure-theoretically, measure-theoretically$^+$, uniformly measure-theoretically$^+$) chaotic then $(X,T)$ is {\rm DC2} (uniformly {\rm DC2}, {\rm DC}{\small$1\!\tfrac12$}, uniformly {\rm DC}{\small$1\!\tfrac12$}) chaotic.
\end{thm}
\begin{proof} 
Let $(\mathscr P_k)$ be a refining \sq\ of partitions such that the diameter of the largest atom in $\mathscr P_k$ decreases to zero with $k$. For each $k$ we define a \sq\ of open sets $U_{k,m}$ ($m\ge 1$) as follows: by regularity of the measure, each atom $P$ of $\mathscr P_k$ can be approximated (in measure) by a \sq\ of its closed subsets, say $(F_{P,m})_{m\ge 1}$. We let 
$$
U_{k,m}=X\setminus \bigcup_{P\in\p_k}F_{P,m}.
$$ 
We have, for every $k$, $\mu(U_{k,m})\underset{m\to\infty}\longrightarrow0$. Also let $s_{k,m}$
denote the (positive) minimal distance between points in different sets $F_{P,m}, F_{P',m}$ with $P,P'\in\mathscr P_k$.

By Theorem \ref{nullset}, if we remove the null set of points which,  for at least one of the sets $U_{k,m}$, do not satisfy the assertion of the ergodic theorem, then in the remaining part the exists an uncountable $(\mathscr P_k)$-scrambled set $E$. We will show that $E$ is {\rm DC2}-scrambled. Let $(x,y)$ be an off-diagonal pair in $E$. For every $\epsilon>0$ the \sq\ of times $n$ when $d(T^nx,T^ny)<\epsilon$ contains the \sq\ of times (of upper density 1) when the points $T^nx,T^ny$ belong to the same atom of $\mathscr P_k$, where $k$ is so large that the diameter of the largest atom of $\mathscr P_k$ is smaller than $\epsilon$. This easily implies that $(x,y)$ satisfies the first requirement for being {\rm DC2}-scrambled. 

Further, there exists $k_0$ and a positive $\eta$ such that $T^nx,T^ny$ belong to different atoms of $\mathscr P_{k_0}$ for $n$'s with upper density at least $\eta$. Let $\delta>0$ be much smaller than $\eta$. Find $m$ so large that the set $U_{k_0,m}$ has measure smaller than $\delta$. If we now remove from the aforementioned \sq\ of times $n$ all the times when at least one of the points $T^nx,T^ny$ belongs to $U_{k_0,m}$, then we are left with a \sq\ of upper density at least $\eta-2\delta$ (still positive) when the two considered points are at least $s_{k_0,m}$ apart. This proves that $(x,y)$ satisfies the second requirement for being {\rm DC2}-scrambled (with the parameters $s=s_{k_0,m}$ and upper density $\eta-2\delta$).
 
If $E$ is uniformly $(\mathscr P_k)$-scrambled, $(\mathscr P_k)^+$-scrambled, or
uniformly $(\mathscr P_k)^+$-scrambled, the same proof yields the corresponding \tl\ scrambling, as in 
the assertion of the theorem.
\end{proof}

There are many examples of {\rm DC2} chaotic systems in which the union of all scrambled sets is a null set for all ergodic measures,
showing that the implication converse to Theorem \ref{topol} need not hold. 
However, if the topological chaos is ``persistent under removing null sets'', it does imply
measure-theoretic chaos, as stated below.

\begin{thm}\label{persistent} Let \xt\ be a \tl\ \ds\ and let $\mu$ be an ergodic \im. Then the system \xsmt\ is measure-theoretically (uniformly measure-theoretically, measure-theoretically$^+$, uniformly measure-\break theoretically$^+$) chaotic if and only if, after removing any set of measure $\mu$ zero, there remains an uncountable {\rm DC2}-scrambled (uniformly {\rm DC2}-scrambled, {\rm DC}{\small$1\!\tfrac12$}-scrambled, uniformly {\rm DC}{\small$1\!\tfrac12$}-scrambled) set.
\end{thm}
\begin{proof}
Necessity follows from the proof of the preceding theorem; the {\rm DC2}-scrambled set (or its variants) has been obtained after removing a specific null set, but we could have additionally removed any other null set as well. We pass to proving sufficiency. Choose a \sq\ $(\mathscr P_k)$ with the diameters of the largest atoms decreasing to zero with $k$, and define the sets $U_{m,k}$ (and the positive numbers $s_{k,m}$) as in the preceding proof. By Theorem \ref{nullset}, it suffices to fix a null set $A$ and find a scrambled set disjoint from $A$. Let $A_0$ be the null set of points which fail the ergodic theorem for at least one of the sets $U_{k,m}$. 
By assumption, there exists an uncountable {\rm DC2}-scrambled set $E$ disjoint from $A\cup A_0$. We will show that $E$ is $(\mathscr P_k)$-scrambled (and it is obviously disjoint from $A$). Take a pair $(x,y)$ of distinct points in $E$ and fix some $k$. Choose an arbitrarily small $\delta>0$ and let $m$ be such that the measure of $U_{k,m}$ is smaller than $\delta$. As we know, the \sq\ of times $n$ when $T^nx$ and $T^ny$ are closer together than $s_{k,m}$ has upper density 1. If we disregard the times when at least one of them falls into $U_{m,k}$, we are left with a \sq\ of upper density at least $1-2\delta$. Note that now, at each of these times, the two points belong to the same atom of $\mathscr P_k$. Because $\delta$ 
is arbitrarily small, we have shown that $T^nx$ and $T^ny$ belong to the same atom of $\mathscr P_k$ for times $n$ of upper density 1. 

We also know that $d(T^nx,T^ny)$ is larger than some positive $s$ for times $n$ with positive upper density. It suffices to pick $k_0$ large enough so that every atom of $\mathscr P_{k_0}$ has diameter smaller than $s$. Then for the same times $n$, $T^nx$ and $T^ny$ must fall into different atoms of $\mathscr P_{k_0}$. This ends the proof for the usual $(\mathscr P_k)$-scrambling. 

The same proof works for the other three variants of chaos.
\end{proof}
 
\begin{rem}\label{totu}
Using Theorem \ref{posset} the following variant of Theorem~\ref{persistent} can be proved: Measure-theoretic chaos (and its respective variants) for an ergodic measure $\mu$ in a \tl\ \ds\ is equivalent to ``$\mu$-ubiquitous {\rm DC2}'' (and its respective variants): an uncountable {\rm DC2}-scrambled (uniformly {\rm DC2}-scrambled, 
{\rm DC}{\small$1\!\tfrac12$}-scrambled, uniformly {\rm DC}{\small$1\!\tfrac12$}-scrambled) set exists within every set of positive measure $\mu$. 
\end{rem}

\begin{rem}\label{dc3}
Here is an analog of {\rm DC3}. Call a pair $(x,x')$ \emph{$(\p_k)^-$-scrambled} if there exists $k_0$ such that
the \sq\ of times $n$ when $T^nx$ belongs to the same atom of $\p_{k_0}$ as $T^nx'$ does not have density (upper and lower densities differ). A system is \emph{measure-theoretically$^-$ chaotic} if, for every refining \sq\ of finite partitions, there exists an uncountable $(\p_k)^-$-scrambled set. Using slight modifications of the proofs presented in this section one can prove that: 1. This notion is an isomorphism-invariant; 2. It suffices to check ``persistence under removing null sets'' and just one refining \sq\ of partitions; 3. It enjoys the ``ubiquitous presence'' property; 4. In \tl\ systems it is equivalent to {\rm DC3} ``persistent under removing null sets''. In particular, as an interesting consequence, we get that {\rm DC3} ``persistent under removing null sets of at least one \im'' is a conjugacy-invariant.
\end{rem}

\section{Measure-theoretic chaos versus entropy}

The most important relation between entropy and chaos is contained in the following theorem, which, combined with our Theorem \ref{topol}
(and the Variational Principle), strengthens several former results (\cite{BGKM}, \cite{D}):

\begin{thm}\label{glowne}
Every ergodic system \xsmt\ with positive Kolmogorov--Sinai entropy is uniformly measure-theoretically$^+$ chaotic.
\end{thm}

\begin{proof} Large part of the proof is identical as in \cite{D}. We now move directly to a certain point of that proof 
(skipping all the arguments that lead to that point). We select a \sq\ of partitions $(\p_k)$ of $X$ in such a way that
$(\p_{2k})$ is a refining \sq\ of partitions (as in the definition of the measure-theoretic chaos), while the odd-numbered
partitions $\p_{2k-1}$ are all equal to one finite partition $\p$ with positive dynamical entropy 
$h_\mu(T,\p)$, which we denote by~$h$. 
We then fix an increasing \sq\ of integers $S=(a_1,b_1,a_2,b_2,a_3,b_3\dots)$. The \sq\ should grow so fast that 
$\frac{b_k}{a_k}$ tends to infinity. We introduce the following notation:

$$
\R_k = \p_k^{[a_k,b_k-1]} := \bigvee_{i=a_k}^{b_k-1}T^{-i}(\p_k),
$$
and 
$$
\R^{\mathsf{odd}}_{1,2k-1}=\bigvee_{i=1}^k\R_{2i-1} \text{ \ \ \ \ and \ \ \ \ }
\mathfrak R = \bigvee_{k=1}^\infty\R_{2k}
$$
(note that $\mathfrak R$ is no longer a finite partition, rather a \emph{measurable partition} which can be identified
with the collection of atoms of the sigma-algebra generated by the partitions involved in the countable join). We will also
denote by $n_k$ the difference $b_k-a_k$. In \cite{D} it is shown that given a decreasing to zero \sq\ of positive numbers $\delta_k$ and a set $X'$ of sufficiently large measure $1-\epsilon_0$ then, if the \sq\ $S$ grows fast enough, there exists an atom $z$ of $\mathfrak R$, a Borel measure $\nu$ supported by $z\cap X'$, a decreasing \sq\ of measurable sets $V_k$, such that

\begin{enumerate}[(A)]
\item if $B$ is an atom of $\R_{1,2k-1}^{\mathsf{odd}}$ contained in $V_k$ then $B$ contains at least $2^{n_{2k+1}(h-\delta_k)}$ atoms of $\R_{1,2k+1}^{\mathsf{odd}}$ contained in $V_{k+1}$ and 
whose conditional measures $\nu_B$ range within $2^{-n_{2k+1}(h\pm\delta_k)}$.\footnote{ 
The measure $\nu$ in \cite{D} is obtained as a disintegration measure $\mu_{yz}$ of $\mu$ with respect to $\Pi\vee \mathfrak R$ 
(where $\Pi$ is the Pinsker sigma-algebra) on an appropriately chosen atom $y\cap z$. The measure is further restricted to the 
intersection of $z$ with the set $X'$. In this paper, this set will be chosen differently than in \cite{D}.  
}
\end{enumerate}
As we shall show in a moment, the statement (A) alone suffices to deduce uniform measure-theoretic$^+$ chaos. 
First of all we remark that if we subtract a null set from $X'$ the statement will still hold (perhaps on a different atom $z$ and for a different \sq\ $S$, but this does not matter).  Thus, in order to complete the proof of Theorem \ref{glowne}, it remains to show that (A) implies the existence of an uncountable uniformly $(\p_k)^+$-scrambled set within the atom $z$.

It is rather easy to see that the statement (A) remains valid if we replace $\nu$ by the conditional measure $\nu_C$, where $C$ is
any set of positive measure $\nu$. We skip the standard argument here  (comp. Fact 1 in \cite{D}). We now represent our space $X$ as a subset of a compact metric space (say, of the unit interval) and using regularity of the measure $\nu$ we can remove a set of small measure from the support of $\nu$ in such a way that all atoms of the partitions $\R_{1,2k-1}^{\mathsf{odd}}$ (for all $k$) intersected with the remaining set $C$ are compact. Replacing $\nu$ by $\nu_C$, we obtain the condition (A) with the additional feature that the atoms $B$ (and those to which $B$ splits) are all compact. This will guarantee that the intersection of any nested chain of such atoms (with growing parameter $k$) is nonempty.\footnote{
This part of the proof -- ensuring nonempty intersections of nested chains -- was handled in \cite{D} differently, by taking closures of the atoms $B$. It could have been handled the same way as we do it here, as well.
}

The remainder of our proof deviates from that in \cite{D}. The main difference is in obtaining separation along a subsequence of upper
density $\eta$ close to 1, (not just positive). This will be achieved not for the partition $\p$ but for $\p^{[0,m-1]}$ with a suitably selected parameter $m$.

At this point we specify the set $X'$. Let $\epsilon_i$ ($i\ge 1$) be a summable \sq\ of positive numbers with small sum 
$\epsilon_0$. Using the Shannon-McMillan Theorem we can find integers $m_i$ and a \sq\ of sets $C_i$, each being a union of less 
than $2^{m_i(h+\epsilon_i)}$ cylinders of length $m_i$, and whose measure exceeds $1-\epsilon_i^2$. Further, using the ergodic 
theorem, we can find $n'_i$ so large that the set of points whose orbits visit $C_i$ more than $n(1-\epsilon_i^2)$ times within 
the first $n$ iterates, for every $n\ge n'_i$, has measure at least $1-\epsilon_i$. For points in a set $X'$ of measure larger 
than $1-\epsilon_0$, this holds for every $i$. 

Fix a number $\eta<1$. Find the smallest parameter $i$ such that, denoting $\epsilon=\epsilon_\eta=\epsilon_i$ and $m=m_\eta=m_i$, 
(the notation $\epsilon_\eta$ and $m_\eta$ will not be used until two pages further) we have $\epsilon<1-\sqrt\eta$ and 
\begin{equation}\label{pipka}
\frac{2H(\sqrt\eta,1-\sqrt\eta)}m+\epsilon(3\#\p+1)<(1-\sqrt\eta)h
\end{equation}
(here $H(p,1-p)$ stands for $-p\log p-(1-p)\log(1-p)$).
Find $k$ such that $b_{2k+1}\ge n'_i$ and, denoting $n=n_{2k+1}$ we have $\frac{2m}n<\epsilon$. Moreover, we require that $2\delta_k$ ($\delta_k$ is the parameter occurring in the condition (A)) and $\frac{\log m}n$ can be added on the left hand side of \eqref{pipka} maintaining the inequality. We remark that the above requirements hold for all sufficiently large $k$. For future reference we let $k_\eta$ be the smallest choice of $k$.

\medskip
Fix an atom $B$ of $\mathscr R^{\mathsf{odd}}_{1,2k-1}$ contained in $V_k$. By (A), this atom contains at least
$2^{n(h-\delta_k)}$ different atoms of $\mathscr R^{\mathsf{odd}}_{1,2k+1}$ contained in $V_{k+1}$. Every such atom has the
form $B\cap A$, where $A$ is an atom of $\mathscr R_{2k+1}$. We will call the atoms $A$ such that $B\cap A$ is nonempty and 
contained in $V_{k+1}$ \emph{good continuations} of $B$. We denote by $\mathcal A(B)$ the collection of good continuations
of $B$ represented as blocks of length $n$, over the alphabet $\p$. We will now count how many blocks $A\in\mathcal A(B)$ may 
disagree with one selected block $A_0\in\mathcal A(B)$ on a smaller than $\eta$ fraction of all subblocks of length $m$.

To do it, we draw the block $A_0$ $m$ times, and we subdivide the $j$th copy ($j=0,1,\dots,m-1$) into subblocks of length $m$ 
by cutting it at positions equal to $j$ mod $m$ (there are at most $\frac nm$ subblocks in each
copy). This diagram shows all subblocks of length $m$ of $A_0$ (plus some incomplete ``prefixes'' and ``suffixes'' at the ends; see Figure \ref{fig1}). 

\begin{figure}
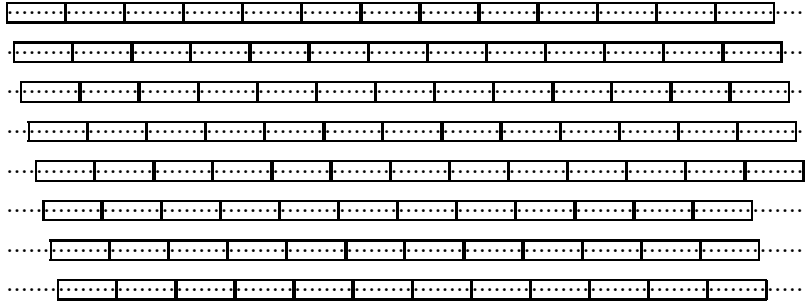

\begin{align*}
&\boxed{\!\!........\!\!}\boxed{\!\!........\!\!}\boxed{\!\!........\!\!}\boxed{\!\!........\!\!}\boxed{\!\!........\!\!}\boxed{\!\!........\!\!}\boxed{\!\!........\!\!}\boxed{\!\!........\!\!}\boxed{\!\!........\!\!}\boxed{\!\!........\!\!}\boxed{\!\!........\!\!}\boxed{\!\!........\!\!}\boxed{\!\!........\!\!}....\\
&.\boxed{\!\!........\!\!}\boxed{\!\!........\!\!}\boxed{\!\!........\!\!}\boxed{\!\!........\!\!}\boxed{\!\!........\!\!}\boxed{\!\!........\!\!}\boxed{\!\!........\!\!}\boxed{\!\!........\!\!}\boxed{\!\!........\!\!}\boxed{\!\!........\!\!}\boxed{\!\!........\!\!}\boxed{\!\!........\!\!}\boxed{\!\!........\!\!}...\\
&..\boxed{\!\!........\!\!}\boxed{\!\!........\!\!}\boxed{\!\!........\!\!}\boxed{\!\!........\!\!}\boxed{\!\!........\!\!}\boxed{\!\!........\!\!}\boxed{\!\!........\!\!}\boxed{\!\!........\!\!}\boxed{\!\!........\!\!}\boxed{\!\!........\!\!}\boxed{\!\!........\!\!}\boxed{\!\!........\!\!}\boxed{\!\!........\!\!}..\\
&...\boxed{\!\!........\!\!}\boxed{\!\!........\!\!}\boxed{\!\!........\!\!}\boxed{\!\!........\!\!}\boxed{\!\!........\!\!}\boxed{\!\!........\!\!}\boxed{\!\!........\!\!}\boxed{\!\!........\!\!}\boxed{\!\!........\!\!}\boxed{\!\!........\!\!}\boxed{\!\!........\!\!}\boxed{\!\!........\!\!}\boxed{\!\!........\!\!}.\\
&....\boxed{\!\!........\!\!}\boxed{\!\!........\!\!}\boxed{\!\!........\!\!}\boxed{\!\!........\!\!}\boxed{\!\!........\!\!}\boxed{\!\!........\!\!}\boxed{\!\!........\!\!}\boxed{\!\!........\!\!}\boxed{\!\!........\!\!}\boxed{\!\!........\!\!}\boxed{\!\!........\!\!}\boxed{\!\!........\!\!}\boxed{\!\!........\!\!}\\
&.....\boxed{\!\!........\!\!}\boxed{\!\!........\!\!}\boxed{\!\!........\!\!}\boxed{\!\!........\!\!}\boxed{\!\!........\!\!}\boxed{\!\!........\!\!}\boxed{\!\!........\!\!}\boxed{\!\!........\!\!}\boxed{\!\!........\!\!}\boxed{\!\!........\!\!}\boxed{\!\!........\!\!}\boxed{\!\!........\!\!}.......\\
&......\boxed{\!\!........\!\!}\boxed{\!\!........\!\!}\boxed{\!\!........\!\!}\boxed{\!\!........\!\!}\boxed{\!\!........\!\!}\boxed{\!\!........\!\!}\boxed{\!\!........\!\!}\boxed{\!\!........\!\!}\boxed{\!\!........\!\!}\boxed{\!\!........\!\!}\boxed{\!\!........\!\!}\boxed{\!\!........\!\!}......\\
&.......\boxed{\!\!........\!\!}\boxed{\!\!........\!\!}\boxed{\!\!........\!\!}\boxed{\!\!........\!\!}\boxed{\!\!........\!\!}\boxed{\!\!........\!\!}\boxed{\!\!........\!\!}\boxed{\!\!........\!\!}\boxed{\!\!........\!\!}\boxed{\!\!........\!\!}\boxed{\!\!........\!\!}\boxed{\!\!........\!\!}.....\\
\end{align*}
		\caption{All subblocks of length $m$ of $A_0$ visualized in $m$ copies of $A_0$.}\label{fig1}
\end{figure}

Imagine another block $A$ treated the same way and suppose it disagrees with $A_0$ on a smaller than $\eta$ fraction of all subblocks. This implies that the fraction of all $m$ copies for which 
$A_0$ and $A$ disagree on a larger than $\sqrt{\eta}$ fraction of subblocks visualized in this copy is at most $\sqrt\eta$. 
In other words, for a fraction of at least $1-\sqrt{\eta}$ of all copies, $A_0$ and $A$ agree on at least a fraction of $1-\sqrt{\eta}$ of the subblocks. Further, we know that, in this diagram, at most $b_{2k+1}\epsilon^2<2n\epsilon$ subblocks represent cylinders not contained in $C_i$. Again, in at least a fraction of $1-\epsilon$ copies the subblocks not contained in $C_i$ constitute a fraction 
smaller than $2\epsilon$. Because $\epsilon+\sqrt\eta<1$, there exists at least one copy where we have both smaller than $2\epsilon$ fraction of subblocks from outside $C_i$ and larger than $1-\sqrt{\eta}$ fraction of agreeing subblocks.
We can now classify all blocks $A$ that we are counting into at most (not necessarily disjoint) 
$$
m \cdot 2^{\frac nm H(\sqrt{\eta},1-\sqrt{\eta})}\cdot 2^{\frac nm H(\epsilon_i,1-\epsilon_i)} \le 2^{\frac nm 2H(\sqrt{\eta},1-\sqrt{\eta})+\log m}
$$
groups depending on the choice of the copy, the choice of necessarily agreeing subblocks (perhaps there will be more), 
and the choice of places reserved for visits in $C_i$ (perhaps not all of them will be used). 

In every group there are at most 
$$
2^{3n\epsilon\#\p}2^{\frac nm\sqrt\eta m(h+\epsilon)} \le 2^{n(\epsilon(3\#\p+1)+h\sqrt\eta)}
$$
blocks. (We allow any symbols from $\p$ on the fraction $2\epsilon$ of subblocks from outside $C_i$ and on the ``prefix''
and ``suffix'' jointly of length not exceeding $2m$, hence constituting another fraction smaller than $\epsilon$. Otherwise,
on a fraction of at most $\sqrt\eta$ of all subblocks we have free choice from the collection of at most $2^{m(h+\epsilon)}$ 
blocks from $C_i$). Multiplying this by the number of groups we get no more than
\begin{equation}
2^{n\bigl(\frac{2H(\sqrt\eta,1-\sqrt\eta)}m+\frac{\log m}n+ \epsilon(3\#\p+1)+h\sqrt\eta\bigr)}< 2^{n(h-2\delta_k)} = \frac{2^{n(h-\delta_k)}}{2^{n\delta_k}}
\end{equation}
blocks.

We have shown that blocks $A$ differing from a selected block $A_0$ on a smaller than $\eta$ fraction of subblocks of 
length $m$ form a negligibly small fraction (at most $2^{-n\delta_k}$) of the family $\mathcal A(B)$. 

\medskip
We are in a position to construct our uncountable uniformly measure-theore\-tically$^+$-scrambled set $E$.
We begin by constructing a family $B_\kappa$, where $\kappa$ ranges over all finite binary words, such that $B_\kappa$
is a (nonempty and closed) atom of $\R^{\mathsf{odd}}_{1,2k+1}$ contained in $V_{k+1}$, where $k$ is the length of $\kappa$ 
(for $k=0$, $\kappa$ is the empty word). We will assure that if $\iota$ extends $\kappa$ to the right then 
$B_{\iota}\subset B_\kappa$. We will also assure an appropriate separation condition. For that we fix a \sq\ $(\eta_k)_{k\ge 1}$
with the following properties: the \sq\ assumes values strictly smaller than but arbitrarily close to 1, each value is
assumed infinitely many times, a value $\eta$ is allowed to occur only for $k\ge k_\eta$. The inductive separation condition
is that if $\kappa$ and $\kappa'$ are binary words of the same length, differing at a position $k_0$, then, for every 
$k\ge k_0$ (up to the length of $\kappa$), the blocks $A_k$ and $A'_k$ differ at a fraction of at least $\eta$ of all
subblocks of length $m_{\eta_k}$, where $A_k$ and $A_k'$ are the blocks appearing at the coordinates $[a_{2i+1},b_{2i+1-1}]$ 
in the symbolic representation of the atoms $B_\kappa$ and $B_{\kappa'}$, respectively. We will do it by induction 
on $k$, in each step we choose two ``children'' of every so far constructed atom of $\R^{\mathsf{odd}}_{1,2k-1}$. 

In step $k=0$ we assign $B_\emptyset$ to be an arbitrarily selected atom of $\R_1$ contained in $V_1$. 
Suppose the task has been completed for some $k-1$, i.e., that we have selected $2^{k-1}$ atoms $B_\kappa$ 
of $\R^{\mathsf{odd}}_{1,2k-1}$, contained in $V_{k-1}$, and pairwise separated as required. We order the $\kappa$'s of length
$k-1$ lexicographically. Take the first atom $B_{\kappa_1}$ (assigned for $\kappa_1=0\dots00$). Choose one good continuation 
$A_0$ of $B_{\kappa_1}$. From every family $\mathcal A(B_\kappa)$ (including $\kappa=\kappa_1$) we eliminate (for future choices) all the atoms $A$ which differ from $A_0$ on a smaller than $\eta_k$ fraction of subblocks of length $m_{\eta_k}$. Since $k\ge k_{\eta_k}$, the preceding estimate applies: every family $\mathcal A(B_\kappa)$ has ``lost'' at most a fraction of $2^{-n_{2k+1}\delta_k}$ of its cardinality. Next we choose $A_1$ from the remaining good continuations of $B_{\kappa_1}$ and again, from each of the families $\mathcal A(B_\kappa)$ we eliminate all the atoms $A$ not sufficiently separated from $A_1$. Again, the losses are negligibly small. We assign $B_{\kappa_10}=B_{\kappa_1}\cap A_0$ and $B_{\kappa_11}=B_{\kappa_1}\cap A_1$ 
(here $\kappa_10$ and $\kappa_11$ denote the two continuations of $\kappa_1$).

Next we abandon $B_{\kappa_1}$ and pass to $B_{\kappa_2}$ ($\kappa_2 = 0\dots01$) and we repeat the procedure choosing two of 
its good continuations, say $A'_0,A'_1$, not eliminated in the preceding steps, each time eliminating for future choices all blocks insufficiently separated from the chosen ones. We proceed until we choose two good continuations for every $\kappa$ of length $k-1$. Note that near the end of this procedure we will have eliminated from each family $\mathcal A(B_\kappa)$ a fraction of at most $2^k\cdot2^{-n_{2k+1}\delta_k}$, which is less than 1 (we decide about the size of $n_{2k+1}$ after fixing $\delta_k$).
Hence the procedure will be possible till the end. This completes the inductive step $k$.

Let now $\kappa$ denote an infinite binary string, while $\kappa_k$ is the prefix of length $k$ of $\kappa$. The atoms 
$B_{\kappa_k}$ form a decreasing \sq\ of compact sets, hence have a nonempty intersection. We select one point from this
intersection and call it $x_\kappa$. The set $E$ is defined as the collection $\{x_\kappa:\kappa\in\{0,1\}^\na\}$.
The following facts are obvious: the set $E$ is uncountable, all its elements belong to the atom $z$ of $\mathfrak R$.
The last fact implies that for each $k$ the orbits of all points from $E$ fall in the same element of $\p_{2k}$ (which 
is a partition in our refining \sq) for all times $n$ belonging to the intervals $[a_{2k'},b_{2k'-1}]$ for all $k'\ge k$.
Because the ratios $\frac{b_{2k}}{a_{2k}}$ tend to infinity, it is clear that such times $n$ have upper density 1 and the
first requirement for $(\p_{2k})^+$-scrambling is verified for all pairs in $E$.

Consider a pair of distinct points from $E$, i.e., $x=x_\kappa$ and $x'=x_{\kappa'}$, where $\kappa\neq\kappa'$. Let $k_0$ denote the
first place where $\kappa$ differs from $\kappa'$. Fix some $\eta<1$ and then, if necessary, replace it by a larger value, so 
that $\eta$ occurs as $\eta_k$ (and then it occurs for infinitely many indices $k$). Pick such a $k$ larger than $k_0$ and
observe the blocks $A_k$ and $A'_k$ representing the atoms of $\mathscr R_{2k+1}$ containing $x$ and $x'$, respectively.
Since $\kappa$ and $\kappa'$ differ at a position smaller than or equal to $k$, the blocks $A_k$ and $A'_k$ have been
selected as either two different continuations of the same atom of $\mathcal R^{\mathsf{odd}}_{1,2k-1}$ or as continuations 
of two different atoms of this partition. In any case, they have been selected one after another in the inductive step $k$, 
which means the latter one (say $A'_k$) was chosen after eliminating all blocks that differ from the former (say $A_k$)
on a smaller than $\eta$ fraction of all subblocks of length $m_{\eta}$. This means that $A_k$ and $A_k'$ differ on
a fraction of at least $\eta$ of such subblocks. Moreover, this is true for infinitely many $k$'s. Because $\frac{b_k}{a_k}$
tends to infinity, this easily implies that the set of times $n$ for which $T^nx$ and $T^nx'$ belong to different atoms
of the partition $\p^{m_\eta}$ has upper density at least $\eta$. This is ``almost'' the second requirement for $(\p_{2k})^+$-scrambling, except that it refers to a wrong partition. In order to replace the partition $\p^{m_\eta}$ by a partition belonging to our refining \sq, i.e., by some $\p_{2k_\eta}$ we apply the same technique as many times before. There exists $k_\eta$ such that $\p_{2k_\eta}$ refines $\p^{m_\eta}$ except on a set of measure $\delta<\!\!<\eta$. If, at the start of the proof, we eliminate points which do not obey the ergodic theorem for the field generated by all the partitions of the form $\p_{2k}$ and $\p^m$, then our points $T^nx$ and $T^nx'$ will belong to different atoms of $\p_{2k_\eta}$ for powers $n$ with upper density at least $\eta-2\delta$ which is as close to 1 as we want. We have shown that the pair $x,x'$ is $(\p_{2k})^+$-scrambled. 
Finally we note is that the assignment $(\eta-2\delta)\mapsto k_\eta$ arising in the proof does not depend on the 
pair $x,x'$ and hence the set $E$ is scrambled uniformly, which concludes the proof of Theorem~\ref{glowne}.
\end{proof}

\begin{rem}\label{rem1}
Let us say that an increasing \sq\ $n_i$ \emph{achieves upper density $\eta$ along a subsequence $N$ of positive integers} if 
$$
\limsup_{N\to\infty}\frac{\#\{i:n_i\le N\}}N \ge\eta.
$$
In the above construction, all upper densities required in the definition of $(\p_k)^+$-scrambling (the upper density 1
of the \sq\ $n_i$ and the upper densities $\eta$, more precisely $\eta-2\delta$, of the \sq s $m_{\eta,i}$) are achieved 
along the sub\sq\ $b_k$ (the right ends of the intervals $[a_k,b_k-1]$). The only constraints on the choice of the \sq\ $S$
(containing $b_k$) concern the speed of its growth, thus $b_k$ could have been selected a sub\sq\ of any a priori given infinite \sq\ of positive integers. We will need this observation in the proof of Theorem \ref{exa}.
\end{rem}

\begin{rem}
We have obtained a specific scrambled set, which, in spite of being ``uniform'' has another property one might call ``synchronic''. 
Let us say that an increasing \sq\ $n_i$ \emph{achieves lower density $\eta$ along a subsequence $N$ of positive integers} if 
$$
\liminf_{N\to\infty}\frac{\#\{i:n_i\le N\}}N \ge\eta.
$$
Clearly, upper density of the \sq\ $a_i$ equals the supremum of all lower densities that the \sq\ achieves along various 
sub\sq s of the positive integers.
For all distinct pairs in our scrambled set $E$ the \sq s of times $n_i$ achieve lower densities 1 along the same sub\sq, 
namely along $b_{2k}$. Further, given $\eta<1$, for all distinct pairs in our scrambled set, the \sq s $m_{\eta,i}$ achieve the lower densities $\eta$ along a common \sq, (for instance, if $\eta$ is assumed as $\eta_k$ then the lower density is achieved along $b_{2k+1}$, where $k$ denotes only these infinitely many integers for which $\eta_k = \eta$). 
\end{rem}

\begin{rem}
In the construction, the \sq s $m_{\eta,i}$ obtained for various values of $\eta$ do not achieve their desired lower densities $\eta$ along the same \sq\ (at least this is not assured). It is so, because with distinct values of $\eta$ we have associated disjoint 
\sq s of indices $k$ such that $\eta_k=\eta$. It is possible to modify the construction to assure the existence of a common 
\sq\ $N$ along which all the \sq s $m_{\eta,i}$ (for varying $\eta$ and varying pairs) would achieve their desired lower 
densities $\eta$ (``synchronic$^+$ scrambling''). This can be done by a different elimination procedure in the construction of the sets $B_\kappa$. In that construction we let the \sq\ $(\eta_k)$ increase to $1$ (without repeating each value infinitely many times) and, in each inductive step (say $k_0)$, we assure appropriate separation (on a fraction $\eta_k$ of all subblocks of length $m_{\eta_k}$) \emph{simultaneously} for all $k\le k_0$. We find this extra property not worth a detailed  proof.
\end{rem}

Combining Theorem \ref{glowne}, and Remark \ref{totu} we obtain the following \tl\ statement, a strengthening of the results from \cite{BGKM} and \cite{D}.

\begin{cor}\label{coro} 
A \tl\ \ds\ \xt\ with positive topological entropy reveals ``ubiquitous chaos {\rm DC}{\small$1\!\tfrac12$}''; an uncountable uniformly 
{\rm DC}{\small$1\!\tfrac12$}-scrambled set exists within every subset of positive measure $\mu$, for every ergodic measure $\mu$ 
with positive Kolmogorov--Sinai entropy.
\end{cor}

For a more complete picture of relations between entropy and our notions of chaos we give an example showing that
Theorem \ref{glowne} cannot be reversed: 

\begin{thm}\label{exa}
There exists a system \xsmt\ with entropy zero and with uniform measure-theoretic$^+$ chaos.
\end{thm}
\begin{proof} As a matter of fact, such an example exists in a paper of Serafin \cite{ser}. Moreover, it is a topological example,
with \tl\ entropy zero, in which we will fix an ergodic measure $\mu$. Because the example was created for different purposes, we 
will need to verify the chaos. This is going to be a tedious task.

\medskip
Let us first say a few words about certain (invertible) systems \xt\ that have an odometer factor. 
Consider a two-row symbolic system, where both rows are bi-infinite \sq s of symbols. The first row contains 
symbols from $\{0,1,2,\dots,\infty\}$, the second row is binary (contains symbols from $\{0,1\}$). The elements $x\in X$ 
obey the following \emph{odometer rule} with respect to an increasing \sq\ $(N_k)$ of integers such that, for each $k$, 
$N_{k+1}$ is a multiple of $N_k$ (this \sq\ is called the \emph{base} of the odometer):
\begin{itemize}
	\item For each $k\ge 1$, the symbols $k'\ge k$ occupy in the first row a periodic set of period $N_k$ having 
	      exactly one element in every period. Such symbols will be called \emph{$k$-markers}. The two-row blocks of length $N_k$
	      starting with a $k$-marker will be called $k$-blocks. Every point in such a system is, for every $k$, a concatenation of the $k$-blocks (see Figure \ref{fig2}).
\end{itemize}
\begin{figure}[ht]
		\includegraphics[width=12truecm]{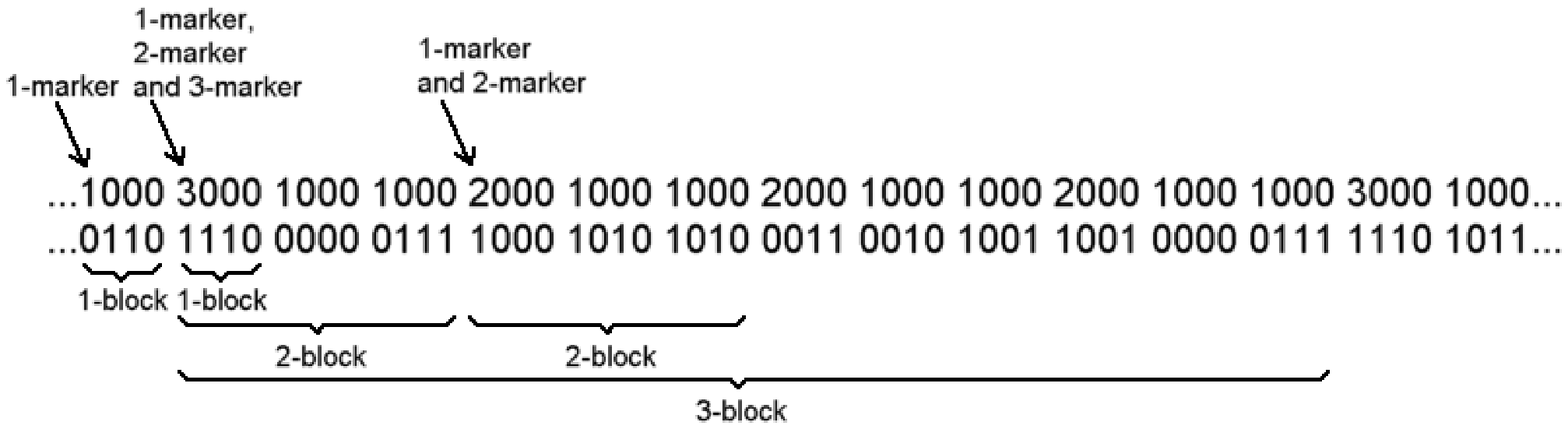}		
		\caption{An element of a system with an odometer factor to base $(N_k)$, with $N_1=4,N_2=12, N_3=48,\dots$.}\label{fig2}
\end{figure}
In such a system we introduce a specific \sq\ of partitions, which we denote by $(\p_k)$ defined in the following manner: 
two points belong to the same atom of $\p_k$ if they have the same and identically positioned \emph{central} $k$-block
(by which we mean the $k$-block covering the coordinate zero). Notice that there are (at most) $N_k\cdot 2^{N_k}$ atoms of $\p_k$ ($N_k$ counts the possible ways the $k$-block is positioned on the horizontal axis, while $2^{N_k}$ is the maximal number of possible ``words'' in the second row of the $k$-block). 

\begin{lem}\label{refine}
The \sq\ of partitions $(\p_k)$ is refining for any \im\ $\mu$. 
\end{lem}
\begin{proof}
By general facts concerning standard probability spaces, it suffices to show that after discarding a null set, the partitions $\p_k$ separate points. It is clear that the partitions separate points $x$ belonging to different fibers of the odometer (i.e., differing in the first row). Notice that the symbol $\infty$ may occur in the first row of an $x\in X$ only one time. Thus, by the Poincar\'e recurrence theorem, the set of elements $x\in X$ in which $\infty$ occurs is a null set for any \im. After discarding this null set, the partitions $\p_k$ also separate points belonging to the same fiber of the odometer (the central $k$-blocks grow with $k$ in both 
directions eventually covering the entire elements). 
\end{proof}

The partitions $\p_k$ have a very specific property (not enjoyed by the ``usual'' partitions of symbolic spaces into blocks occurring at fixed positions): if $x,x'$ belong to the same atom of $\p_k$ then $T^ix,T^ix'$ belong to the same atom of $\p_k$, for $i$ ranging in an interval of integers of length $N_k$ containing $0$ (namely as long as shifting by $i$ positions maintains the coordinate zero within the same $k$-block). In particular, if $k'>k$ and $x,x'$ have identically positioned 
$k'$-blocks, say $x[a,b]$ and $x'[a,b]$ are $k'$-blocks, then the percentage of times $i\in [a,b]$ when $T^ix$ and $T^ix'$ belong to the same (different) elements of $\p_k$ equals the percentage of agreeing (disagreeing) component $k$-blocks in the $k'$-blocks $x[a,b]$ and $x'[a,b]$. We will refer to this property at the end of the proof.

Notice that if two points belong to different fibers of the odometer factor then their orbits are separated by some
$\p_k$ at all times. Thus every $(\p_k)$-scrambled set (if one exists) is contained in one fiber of the odometer. In particular, the uniform measure-theoretic$^+$ chaos is equivalent to the condition that $\nu$-almost every element 
$y$ of the odometer, where $\nu$ is the unique \im\ on the odometer, has the property that after discarding any null 
set for the disintegration measure $\mu_y$, the fiber of $y$ contains an uncountable uniformly $(\p_k)^+$-scrambled set.
\medskip

We will now reproduce the construction of the example from \cite{ser}. 
At first, we introduce some notation. If $\mathcal B\subset \{0,1\}^k$ (i.e., $\mathcal B$ is a family of selected binary blocks of length $k$) and $q\in\na$ then by $\mathcal B^q$ we mean the family of all concatenations of $q$ elements from $\mathcal B$, $\mathcal B^q=\{B_1B_2\dots B_q: B_1,B_2,\dots,B_q\in\mathcal B\}$ and by $\mathcal B_\mathsf{rep}$ we will mean $\{BB: B\in\mathcal B\}$ (the family of all two-repetitions of elements of $\mathcal B$). Notice that $\mathcal B_\mathsf{rep}\subset\mathcal B^2$, $\#\mathcal B^q = (\#\mathcal B)^q$ and $\#\mathcal B_\mathsf{rep} = \#\mathcal B$.

We fix a \sq\ of integers $(q_k)_{k\ge 1}$, larger than 1 (for the purposes of this example it suffices to take $q_k=2$ for all $k$). The products $q_1q_2\cdots q_k$ will be denoted by $p_k$. Next we define inductively two families of binary blocks: 
$\mathcal B_1 = \{0,1\}^{q_1}$, $\mathcal C_1 = (\mathcal B_1)_{\mathsf{rep}}$, and, for $k\ge 2$, $\mathcal B_k = (\mathcal C_{k-1})^{q_k}$, $\mathcal C_k = (\mathcal B_k)_{\mathsf{rep}}$. The blocks in $\mathcal B_k$ have length $p_k2^{k-1}$ and those in $\mathcal C_k$ have length $N_k= p_k2^k$. We have $\#\mathcal B_k = \#\mathcal C_k = 2^{p_k}$. 

Now we can define the system \xt\ as a two-row symbolic system with the odometer $G$ to base $(N_k)$ encoded 
in the first row and such that the ``word'' in the second row of every $k$-block is a block from $\mathcal C_k$ (see Figure \ref{fig3}.)

\begin{figure}[ht]
\begin{align*}
\dots &100\,000\ 200\,000\ 100\,000\ 100\,000\ 100\,000\ 100\,000\ 100\,000\ 200\,000\ 100\,000\dots\\
\dots &101\,101\ \underset{\ \ \ \in\mathcal C_2}{\underbrace{\underset{\ \ \ \in\mathcal C_1}{\underbrace{010\,010}}\ 
\underset{\ \ \ \in\mathcal C_1}{\underbrace{110\,110}}\ \underset{\ \ \ \in\mathcal C_1}{\underbrace{000\,000}}\ \underset{\ \text{repetition}}{\underbrace{010\,010\ 110\,110\ 000\,000}}}}\ 111\,111\ 011\,011\dots
\end{align*}
		\caption{An element of our system (here $q_1=q_2=3$ hence $N_1=6$, $N_2=36$).}\label{fig3}
\end{figure}

The system so defined has \tl\ entropy zero: the first row system has \tl\ entropy zero because it is an odometer and the second row factor has entropy zero by an easy counting argument: the logarithm of the cardinality of the family of $k$-blocks (equal to $p_k$) grows much slower than their length (equal to $p_k2^k$). (The entire system has entropy zero, as a \tl\ joining of two systems with entropy zero.)

In order to define an \im\ $\mu$ it suffices to declare all atoms of the partition $\p_k$ (corresponding to the $k$-blocks positioned around coordinate zero) to have equal measures. (The measure of an atom then equals $\frac1{\#\p_k}=\frac1{N_k2^{p_k}}$.) We skip the standard proof that this indeed determines a shift-\im\ which is ergodic.

Let $H$ be the odometer to base $(p_k)$ represented similarly as $G$, in form of a symbolic system over the alphabet $\{0,1,\dots,\infty\}$. Notice that since, for each $k$, $N_k$ is a multiple of $p_k$, $H$ is a \tl\ factor of $G$; the factor map (which we denote by $\phi$) consists in inserting more $k$-markers (we skip the obvious details). The unique \im\ $\nu$ supported by $G$ is sent
to the unique \im\ $\xi$ supported by $H$.

For $\nu$-almost every $y\in G$ we will now describe a measurable bijection $\pi_y$ between the fiber of $y$ in our system $X$ and the fiber of $\phi(y)$ in the  direct product $(H\times\{0,1\}^\z)$, sending the measure $\mu_y$ to the Bernoulli measure $\lambda=\{\frac12,\frac12\}^\z$. 
The rigorous definition of the map and the proof of the correspondence of measures, although completely elementary, are lengthy and 
not very interesting. Instead we provide a slightly informal description.

Suppose we want to create a $k$-block appearing in the system $X$. Since the first row is determined (up to the value of the ``leading marker'' which can be any number $k'\ge k$) we only need to write a block belonging to $\mathcal C_k$. Suppose we write it from 
left to right. While doing this, we encounter two kinds of positions: those which can be filled completely arbitrarily, independently of what was filled earlier (we call them the ``free positions''), and other, where we have ``forced repetitions''  of something that was filled earlier. There are exactly $p_k$ free positions and each of them is repeated $2^k-1$ times (jointly determining $2^k$ positions, see Figure \ref{fig4}). 

\begin{figure}[ht]
\begin{align*}
&3000\,1000\,1000\,1000\,2000\,1000\,1000\,1000\,2000\,1000\,1000\,1000\,2000\,1000\,1000\,1000\\
&\underline{01}01\,\underline{11}11\,0101\,1111\,\underline{0\overline0}0\overline0\,\underline{10}10\,0\overline00\overline0\,1010\,0101\,1111\,0101\,1111\,0\overline00\overline0\,1010\,0\overline00\overline0\,1010
\end{align*}
		\caption{An example of a 3-block (for $q_1=q_2=q_3=2$). The ``free positions'' are underlined and the positions determined 
		by the free position 18 are overlined.}\label{fig4}
\end{figure}

\noindent Given a $k$-block $B$ appearing in $X$, the free positions read from left to right constitute a new two-row block of length $p_k$, with markers positioned as in a $k$-block over the odometer $H$. We have just described a bijection (denoted $\pi_k$) between all $k$-blocks of $X$ and all $k$-blocks appearing in $(H\times\{0,1\}^\z)$. 

We need to take a closer look at the map $\pi_k$. Every $(k\!+\!1)$-block $C$ is a concatenation of $2q_k$ $k$-blocks (half of which
is repeated): 
$$
C=B^{(1)}B^{(2)}\dots B^{(q_k)}B^{(1)}B^{(2)}\dots B^{(q_k)}.
$$ 
The reader will easily verify that then
$$
\pi_{k+1}(C) = \pi_k(B^{(1)})\pi_k(B^{(2)})\dots\pi_k(B^{(q_k)}).
$$

Fix some $y\in G$ and let $x$ belong to the fiber of $y$ (i.e., $x$ has $y$ in the first row). For each $k$, let $B_k(x)$ denote 
the central $k$-block in $x$. Find the coordinate zero in the block $B_{k+1}(x)$. Since $B_{k+1}(x)=B^{(1)}B^{(2)}\dots B^{(q_k)}B^{(1)}B^{(2)}\dots B^{(q_k)},$ there exists an index $i$ such that the coordinate zero falls within (one of two copies) of the $k$-block $B^{(i)}$ (that is to say, $B_k(x)=B^{(i)}$). In such case we place the image block $\pi_{k+1}(B_{k+1}(x))$ (which equals $\pi_k(B^{(1)})\pi_k(B^{(2)})\dots\pi_k(B^{(q_k)})$) along the horizontal axis in such a way that the zero coordinate falls within 
the subblock $\pi_k(B^{(i)})$. The precise location of the coordinate zero within the latter block is established analogously, by repeating the same procedure for $k,k-1, k-2$, etc. In this manner we obtain a consistent family of positioned blocks $\pi_k(B_k(x))$ growing either in one or in both directions around the coordinate zero (and this behavior depends only on the first row $y$ of $x$). 
By an easy argument, $\nu$-almost surely, the blocks $\pi_k(B_k(x))$ grow around the coordinate zero in both directions, eventually determining a two-sided \sq, which we denote by $\pi_y(x)$. This completes the definition of the map $\pi_y$. We skip the verification
of the fairly obvious fact that for $y$'s for which the map $\pi_y$ is defined, $\pi_y$ is a bijection between the fiber of $y$ and $\{\phi(y)\}\times\{0,1\}^\z$, sending $\mu_y$ to $\lambda$. 

For $y$'s for which the map $\pi_y$ is well defined, the partition $\p_k$ (determined by the position and contents of the central $k$-block) restricted to the fiber of $y$ maps by $\pi_y$ to an analogous partition of $H\times\{0,1\}^\z$ (which we denote by $\Q_k$),
restricted to $\{\phi(y)\}\times\{0,1\}^\z$. Since the product system $H\times\{0,1\}^\z$ with the product measure $\xi\times\lambda$ has obviously positive entropy, by Theorem \ref{glowne}, it is uniformly measure-theoretically$^+$ chaotic. As we have already mentioned, this means that $\xi$-almost every $z\in H$ (we let $H'$ be the corresponding set of full measure) has the property that after removing any $\lambda$-null set, $\{z\}\times\{0,1\}^\z$ contains an uncountable uniformly $(\Q_k)^+$-scrambled set. We would like to deduce that this property passes to $X$. We will use the map $\phi$ between the odometers and the maps $\pi_y$ on the fibers. For instance we immediately define $G'=\phi^{-1}(H')$ and we note that $G'$ is a set of full measure $\nu$. The difficulty is that the combined map (from $X$ to $H\times\{0,1\}^\z$) is not shift-equivariant (it cannot be, because a system with positive entropy cannot be a factor of a system with entropy zero), hence the ``percentage of agreement/disagreement'' along two orbits is not automatically preserved. We must check it ``manually''. 

Notice that, since every free position in a $k$-block of $X$ has the same number of ``forced repetitions'', the percentage of entries where two such $k$-blocks differ (equal) is the same as the percentage of entries where their images by $\pi_k$ differ (equal). 
Moreover, this ``preservation of percentage'' passes to higher blocks: Pick some $k'>k$. Every $k'$-block $C$ is a concatenation of 
$k$-blocks; some of them are ``free'', and some are ``forced repetitions''. Each of the ``free'' component $k$-blocks has the same 
number of ``forced repetitions''. The reader will easily observe that the image $\pi_{k'}(C)$ equals the concatenation of the 
images by $\pi_k$ of the ``free'' component $k$-blocks. Thus, for two $k'$-blocks, say $C,D$, the percentage of the component 
$k$-blocks which are the same (different) in $C$ and $D$ is the same as the percentage of the 
component images of the $k$-blocks which are the same (different) in $\pi_{k'}(C)$ and $\pi_{k'}(D)$. 

By an easy argument involving the Borel-Cantelli Lemma, $\nu$-almost every $y$ has the following property: there exists a sub\sq\ 
$k_l$ (depending on $y$) such that the relative position of the zero coordinate within the central $k_l$-block divided by its 
length $N_{k_l}$ converges with $l$ to zero. Similarly, for $\nu$-almost every $y$, the same holds (and we can assume, along 
the same \sq\ $k_l$) for the images $\pi_y(x)$. 

We now fix some $y\in G'$ for which the map $\pi_y$ is defined and which fulfills the above two ``almost sure'' conditions. Let 
$A$ be a $\mu_y$-null set and let $A'=\pi_y(A)$ (which is a null set for the Bernoulli measure). We already know that $\{\phi(y)\}\times(\{0,1\}^\z\setminus A')$ contains an uncountable uniformly $(\Q_k)^+$-scrambled set $E'_A$. Moreover, by Remark \ref{rem1} we can arrange such a scrambled set $E'_A$ that all the upper densities required in the definition of 
scrambling are achieved along the sub\sq\ $p_{k_l}$. 

Let $E_A$ be the preimage by $\pi_y$ of $E'_A$. Clearly, $E_A$ is uncountable and disjoint from $A$. It remains to show is that 
$E_A$ is uniformly $(\p_k)^+$-scrambled. This, however, is an almost immediate consequence of the following facts (we leave the easy 
deduction to the reader):
\begin{itemize}
	\item the upper densities required for scrambling of $E'_A$ are all achieved along the sub\sq\ $p_{k_l}$, 
	\item for large $l$ the central $k_l$-blocks in $x,x'$ (in the fiber of our selected $y$) start ``nearly'' at the coordinate zero,            which implies that the percentage of times $i\in [0,N_{k_l}-1]$ when $T^ix$ and $T^ix'$ belong to the same (different)              atoms of $\p_k$ is nearly the same as for $i\in [a,b]$, where $a,b$ are the ends (common for both points) of the central
	      $k_l$-block. An analogous statement holds for $\pi_y(x)$ and $\pi_y(x')$ and the times $i\in[0,p_{k_l}-1]$ (of course, the              ends $a,b$ of the central $k_l$-blocks are now different);
  \item due to the specific property of the partitions $\p_k$, the percentage of times $i\in[a,b]$, when $T^ix$ and $T^ix'$ belong          to the same (different) atoms of $\p_k$ equals the percentage of agreeing (disagreeing) $k$-blocks in the central $k_l$-blocks          of these points. An analogous statement holds for the points $\pi_y(x),\pi_y(x')$ and the partitions $\Q_k$; 
	\item for $k'>k$ the percentage of agreeing (disagreeing) $k$-blocks within the central $k_l$-blocks of $x,x'$ equals the percentage          of agreeing (disagreeing) $k$-blocks within the central $k_l$-block of $\pi_y(x), \pi_y(x')$. \qedhere
\end{itemize}
\end{proof}

\section{Open problems}

\begin{que}{\bf 1.}
As we have already mentioned, we do not know whether {\rm DC2} (or uniform {\rm DC2}) persistent under removing null sets implies {\rm DC}{\small$1\!\tfrac12$}.
Similarly, we have no examples showing that measure-theoretic$^+$ chaos is essentially stronger than 
the measure-theoretic chaos. 
\end{que}

\medskip
As the results of \cite{DL1} and \cite{DL2} show, some properties necessary for positive \tl\ entropy (such as the existence of asymptotic or forward mean proximal pairs\footnote{A pair $x,x'$ is \emph{forward mean proximal} if there exists a sequence $n_i$ of density~1 along which $d(T^{n_i}x,T^{n_i}x')\to 0$ (in other words $\Phi_{x,y}(0)=1$). This terminology goes back to Ornstein and Weiss \cite{OW} (or even to an earlier work of Furstenberg). It is regretful that such pairs are not called ``mean asymptotic''. ``Mean proximal'' fits much better to pairs for which the \sq\ $n_i$ has \emph{upper} density 1 (i.e., with $\Phi^*_{x,y}(0)=1$), which is the first condition in {\rm DC1} and {\rm DC2} scrambling. With the present terminology, we have no good name for such pairs.
}), 
if inherited by all \tl\ extensions of the system, become also sufficient. Thus it seems reasonable to ask the following

\begin{que}{\bf 2.}
Is it true that a \tl\ \ds\ whose every \tl\ extension is {\rm DC2} ({\rm DC}{\small$1\!\tfrac12$}) chaotic has positive \tl\ entropy? Is it true that every ergodic system whose every (measure-theoretic) extension is measure-theoretically (measure-theoretically$^+$) chaotic has positive Kolmogorov--Sinai entropy?
\end{que}

We were unable not only to resolve the above, but even to disprove the following: it might happen that unlike \tl\ chaoses, the measure-theoretic chaos passes to extensions. So we have another question, in a sense opposite to the preceding one:

\begin{que}{\bf 3.}
Suppose \xsmt\ is a (measure-theoretic) factor of \ycns\ and that the former system reveals measure-theoretical chaos (in one 
of the four discussed versions). Does that imply the same chaos for the latter system? Similarly, if \xt\ is a \tl\ factor
of \ys\ and \xt\ reveals DC2 persistent under removing null sets. Does that imply DC2 for \ys?
\end{que}

\bigskip\noindent
Institute of Mathematics and Computer Science, Wroclaw University of Technology, Wybrze\.ze Wyspia\'nskiego 27, 50-370 Wroc\l aw, Poland

\noindent
downar@pwr.wroc.pl

\smallskip\noindent
Institut des Sciences de l'Ing\'enieur de Toulon et du Var, 
Avenue G. Pompidou, B.P. 56, 83162 La Valette du Var Cedex, France

\noindent
yves.lacroix@univ-tln.fr

\end{document}